\documentclass{article}

\begin{document}

\title{Volume estimates for K\"ahler-Einstein metrics: the three dimensional case}
\author{X-X. Chen and S. K. Donaldson}
\date{\today}
\maketitle

\newcommand{\tg}{\tilde{g}}
\newcommand{\cA}{{\cal A}}
\newcommand{\bQ}{{\bf Q}}
\newcommand{\cX}{{\cal X}}
\newcommand{\bC}{\mbox{${\bf C}$}}
\newcommand{\bR}{\mbox{${\bf R}$}}
\newcommand{\Var}{{\rm Var}}
\newcommand{\Av}{{\rm Av}}
\newcommand{\Vol}{{\rm Vol}}
\newcommand{\Dim}{{\rm Dim}}
\newcommand{\cO}{{\cal O}}
\newcommand{\cW}{{\cal W}}
\newcommand{\cL}{{\cal L}}
\newcommand{\Tr}{{\rm Tr}}
\newcommand{\Zmax}{Z_{{\rm max}}}
\newcommand{\Zmin}{Z_{{\rm min}}}
\newcommand{\Ch}{{\rm Ch}}
\newcommand{\bP}{\mbox{${\bf P}$}}
\newcommand{\uA}{\mbox{${\underline{A}}$}}
\newcommand{\uM}{\mbox{${\underline{M}}$}}
\newcommand{\um}{\mbox{${\underline{m}}$}}
\newcommand{\ur}{\mbox{${\underline{r}}$}}
\newtheorem{Goal}{Goal}
\newtheorem{thm}{Theorem}
\newtheorem{prop}{Proposition}
\newtheorem{lem}{Lemma}
\newtheorem{defn}{Definition}
\newcommand{\oK}{\overline{K}}
\newcommand{\dbd}{\sqrt{-1} \partial\overline{\partial}}
\newcommand{\ulambda}{\underline{\lambda}}
\newcommand{\olambda}{\overline{\lambda}}
\newcommand{\Riem}{{\rm Riem}}
\newcommand{\Ric}{{\rm Ric}}


\section{Introduction}

This is the first of a series of papers in which we obtain estimates for the volume of certain  subsets of K\"ahler-Einstein manifolds. These estimates form the main analytical input for an approach to  general existence
questions \cite{kn:D1}. Let $(M,g)$ be any compact Riemannian manifold and $r>0$. We let $K_{r}\subset M$ be the set of points $x$ where $\vert \Riem\vert\geq r^{-2}$ and write $Z_{r}$ for the $r$-neighbourhood of $K_{r}$. Thus any point of the complement $M\setminus Z_{r}$ is the centre of a metric $r$-ball on which the curvature is bounded by $r^{-2}$. Re-scaling lengths by a factor $r^{-1}$, this becomes, in the re-scaled metric, a unit ball on which the curvature is bounded by $1$. If in addition (as will be the case in our situation),
we have control of the local injectivity radius we can say that, at the length scale $r$, the complement $M\setminus Z_{r}$ consists of \lq\lq good points'' with neighbourhoods of bounded geometry. Our aim is to derive estimates for the {\it volume} of the \lq \lq bad'' set $Z_{r}$: the set where the geometry need not be standard at this scale.

Let $g$ be a K\"ahler-Einstein metric so $\Ric=\lambda g $ for constant $\lambda$. Then one can find in each complex dimension $n$ constants $a_{n}, b_{n}$, depending only on $n$, so that pointwise on the manifold,
\begin{equation}  \vert \Riem \vert^{2} d\mu = (a_{n} c_{1}(\Riem)^{2} + b_{n} c_{2}(\Riem))\wedge \omega^{n-2} , \end{equation}
where we write $c_{i}(\Riem)$ for the standard integrand defining the Chern classes in Chern-Weil theory. (These constants do not depend on $\lambda$. As usual we write $\omega$ for the 2-form corresponding to $g$.)
Thus on a compact manifold $M$
\begin{equation} \int_{M} \vert \Riem \vert^{2} = \langle a_{n} c_{1}^{2}\cup \omega^{n-2}+ b_{n} c_{2}\cup \omega^{n-2}, M \rangle . \end{equation}
The right hand side is a topological invariant, determined by the Chern classes of $M$ and the K\"ahler class,  which we will denote for brevity by $E(M)$. This identity shows that the curvature of the K\"ahler-Einstein metric cannot be very large on a  set of large volume in $M$. More precisely we have an obvious estimate for the volume of $K_{r}$

$$  \Vol (K_{r}) \leq E(M) r^{4}. $$

The goal of our work  is to estimate the volume of $Z_{r}$ rather than $K_{r}$. We consider K\"ahler-Einstein metrics  with non-negative Ricci curvature and in this first paper we restrict to complex dimension $3$.   We consider a compact K\"ahler-Einstein  3-fold  $(M,g)$   with  $\Ric =\lambda g\ , \lambda \geq 0$.  We suppose that the metric satisfies the  condition that 
\begin{equation}  \Vol M \geq \kappa \ {\rm Diam}(M)^{6} \end{equation}
for some $\kappa>0$. The Bishop-Gromov comparison theorem implies that for all metric balls $B(x,r)$ with $r\leq {\rm Diam}(M)$,
\begin{equation}    {\rm Vol}\ B(x,r) \geq \kappa r^{6}. \end{equation}

Our main result is the following.
\begin{thm}

In this situation $$  \Vol(Z_{r})\leq C \left( E(M) r^{4}+ b_{2}(M) r^{6}\right) $$
where $C$ depends only on $\kappa$ and $b_{2}(M)$ is the second Betti number of $M$. 
\end{thm}

Notice that the statement here is scale invariant (as of course it has to be).
When $c_{1}>0$ (which is the case we have mainly in mind) the bound (3) follows from Myers' Theorem,  with  a constant $\kappa$ determined by  topological data.
 If we fix the scale by requiring that $\Ric=g$ then the content of the theorem is the bound
\begin{equation}   \Vol(Z(r))\leq C r^{4}. \end{equation}

We will also establish a small extension, which is probably not optimal.
\begin{thm}
With notation as in Theorem 1, there is a constants $C' $ such that for all $r$ there is a connected open subset $\Omega'\subset M\setminus Z_{r}$ with $$\Vol(M\setminus \Omega')\leq C (E(M) r^{18/5}+ b_{2}(M) r^{6})$$
\end{thm}
(In applications such as in \cite{kn:D1} any bound $O(r^{\mu})$ with $\mu>2$ will suffice.)

We can  consider the same questions for K\"ahler-Einstein metrics
 in any complex dimension $n$.
When $n=2$ the proofs are much easier because of the {\it scale invariance} of the $L^{2}$ norm of the curvature in that dimension. In general for a ball $B(x,r)$ (with $r\leq {\rm Diam}(M)$) we define the \lq\lq normalised energy''
\begin{equation}    E(x,r)= r^{4-2n} \int_{B(x,r)} \vert \Riem \vert^{2} d\mu, \end{equation}
which is  scale invariant. (When we want to emphasise the dependence on the metric we write $E(x,r,g)$.) Notice that in our situation we have
$$  \kappa r^{2n} \leq \Vol (B(x,r)) \leq \omega_{2n} r^{2n} $$
where $\omega_{2n}$ is the volume of the unit ball in $\bR^{2n}$.  So it is essentially the same to normalise by the appropriate power of the volume of the ball.   Normalised energy functionals of this kind appear in many other contexts in differential geometry, for example the theories of  harmonic maps and Yang-Mills fields. In these two theories a crucial {\it monotonicity} property holds. This is the statement that with a fixed centre the normalised energy is a decreasing function of $r$. If this monotonicity property held in our situation for K\"ahler-Einstein metrics, the proof of our theorem would be relatively straightforward. (Of course when $n=2$ the monotonicity is obviously  true.) The main work in this paper is to establish a result which can be seen as an \lq\lq approximate monotonicity'' property. To state this cleanly, let us say that an open subset  $U\subset M$ {\it carries homology}
if the inclusion map $H_{2}(\partial U,\bR)\rightarrow H_{2}( U,\bR)$ is not surjective. Note that if $V\subset U$  carries homology then so does $U$ and that if $U_{1}, U_{2}\dots, U_{p}\subset M$ are disjoint domains which each carry homology then $p$ cannot exceed the second Betti number of $M$.  

\begin{thm} With the same hypotheses as Theorem 1, for each $\epsilon>0$ there is a $\delta>0$ such that if $B(x,r)\subset M$ is a metric ball ($r\leq {\rm Diam}(M)$) which does not carry homology and $E(x,r)\leq \delta$ then 
for any $y\in B(x,r/2)$ and $r'\leq r/2$ we have $E(y,r')\leq \epsilon$.
\end{thm}
The function $\delta(\epsilon)$ depends only on the non-collapsing constant $\kappa$.  

In Section 2 of this paper we prove Theorem 3. The proof is an application of the extensive theory, due to Anderson, Cheeger, Colding and Tian, of Gromov-Hausdorff limits  of Riemannian manifolds with  lower bounds on Ricci curvature. In particular we make use of  deep results of Cheeger, Colding and Tian on
codimension 4 singularities and of Cheeger and Colding on tangent cones at infinity.   Given Theorem 3, the deduction of Theorem 1, which we do  in Section 3, is fairly straightforward. We will first prove a \lq\lq small energy'' result, as follows. 
\begin{prop} With the same hypotheses as Theorem 1, there are $\delta_{0}, K$ such that if $B(x,r)\subset M$  (with $r\leq {\rm Diam}(M)$) is a ball which does not carry homology and $E(x,r)\leq \delta_{0}$ then $\vert \Riem \vert \leq K r^{-2}$ on the interior ball $B(x,r/3)$. 
\end{prop}

Theorem 1 follows from a straightforward covering argument. In Section 4 we conclude with some remarks and discussion.

In the sequel to this paper we will extend the  results to all dimensions using a rather different argument, making more use of the complex structure and developing ideas of Tian in \cite{kn:T1}. (Tian has informed us that, using these ideas, he obtained related results some while ago.) This argument also gives another approach to the three-dimensional case here. But it appears to us worthwhile to write down both proofs. 

The authors have had this paper in draft form since early 2010. Recently Cheeger and Naber have posted a preprint \cite{kn:CN} establishing these volume estimates in all dimensions, and including results in a more general Riemannian geometry setting. Their approach is somewhat different and it seems valuable to have this variety of arguments in the literature.

\section{Proof of Theorem 2}

\subsection{Cones with small energy}

One foundation of our proof of Theorem 1  is a result of Cheeger, Colding and Tian which states, roughly speaking,  that the formation of codimension-4 singularities requires a definite amount of energy.   Let $G\subset U(2)$ be a finite group acting freely on $S^{3}$ and $n\geq 2$. Consider the unit ball in the metric product $\bC^{n-2}\times \left( \bC^{2}/G\right)$ centred at $(0,0)$. Let $V$ be a K\"ahler manifold of complex dimension $n$ with nonnegative Ricci curvature  and $B'$ be a unit ball centred at $p\in V$. Write  $d_{GH}(B,B')$ for the based Gromov-Hausdorff distance  
\begin{prop}(\cite{kn:CCT}, Theorem 8.1) There are $\alpha_{n}, \eta_{n}>0$ such that
if $d_{GH}(B,B')\leq \alpha_{n}$ and $G$ is non-trivial then
$$  \int_{B'} \vert {\rm Riem} \vert^{2} \geq \eta_{n}. $$
\end{prop}

Using this we now prove a result about metric cones in complex dimension $3$ with small energy.
\begin{prop}Let $Y$ be a metric cone with vertex $O$. Suppose that the based space $(Y,O)$ is the Gromov-Hausdorff limit of based K\"ahler-Einstein manifolds, (of complex dimension $3$)
$(X_{n}, g_{n}, O_{n})$ with $\Ric(g_{n})=\lambda_{n} g_{n}$ where $\lambda_{n}\geq 0$ and $\lambda_{n}\rightarrow 0$. Suppose that
the $X_{n}$ satisfy a non-collapsing condition (4), with fixed $\kappa>0$. Then there is a $\theta>0$ such that if  $E(O_{n}, 2, g_{n})\leq \theta $ for all $n$ then $Y$ is smooth away from the vertex.
\end{prop} 

To see this we argue by contradiction. Suppose there is a sequence of such examples $Y_{m}$, with fixed $\kappa$ and with $\theta_{m}\rightarrow 0$. Taking a subsequence, we can suppose these have a based Gromov-Hausdorff limit $Y_{\infty}$ and by a diagonal argument this is the limit of a sequence of smooth based manifolds $X_{n}$, as above, with $E(O_{n}, 2, g_{n})\rightarrow 0$. The non-collapsing condition means that, according to the Cheeger, Colding, Tian theory, $Y_{\infty}$ is a smooth Ricci-flat K\"ahler manifold, of real dimension 6, outside a set ${\cal S}\subset Y_{\infty}$ of Hausdorff dimension at most $2$. Further, because of the complex structures present, if the dimension of ${\cal S}$ is strictly less than $2$ it must be $0$ and this means that the singular set contains at most the vertex $O$, since it is invariant under the dilation action on the cone. So it suffices to show that $\dim {\cal S}<2$. If, on the contrary, the dimension is $2$ then, again by the general theory, there is some point with a tangent cone of the form $ \bC \times \bC^{2}/G$ for a nontrivial $G$. (In fact this is true for almost all points of the singular set, with respect to Hausdorff 2-measure.) By the invariance under the dilation action we can suppose this point is at distance $1$, say, from the vertex. By the definition of tangent cone and of Gromov-Hausdorff convergence we can find a fixed small $\rho$ such that for all large $n$ there is a point $p_{n}\in X_{n}$ such that the distance in $X_{n}$ from $p_{n}$ to $O_{n}$ is approximately $1$ and the  Gromov-Hausdorff distance between the $\rho$-ball in $X_{n}$ centred at $p_{n}$ and the $\rho$-ball in the model space $\bC\times \bC^{2}/G$ is less than $\alpha_{3}\rho$. But, rescaling the result above, this implies that the integral of $\vert \Riem\vert^{2}$ over this $\rho$-ball in $X_{n}$ is at least $\eta_{3}\rho^{2}$ which gives a contradiction.

\begin{prop}
With the same hypotheses as in Proposition 2, for any $\sigma , L$ we can choose $\theta $ so that for all large $n$ and points $x$ in $X_{n}$ with $L^{-1}\leq d(x,O_{n})\leq L $ we have
$$ \vert \Riem \vert\leq \sigma d(x,O_{n})^{-2}. $$
\end{prop}
 This follows from the same argument as above, and general theory.(By results of Anderson \cite{kn:MA}, at points where the Gromov-Hausdorff limit is smooth the metrics converge in $C^{\infty}$.)

\begin{prop}
With the same hypotheses as in Proposition 2, we can choose $\theta$ so that the cone $Y$ must have the form $\bC^{3}/\Gamma$ for some $\Gamma\subset U(3)$ acting freely on $S^{5}$.
\end{prop}

We know that $Y$ has a K\"ahler metric with zero Ricci curvature, so this follows from the {\it rigidity} of $\bC^{3}/\Gamma$, among such cones. To give a direct argument we use the fact that the curvature tensor of a K\"ahler-Einstein metric satisfies an identity of the form
\begin{equation}   \nabla^{*}\nabla \Riem = \Riem *\Riem,\end{equation}
for a certain bilinear algebraic expression $*$. With suitable normalisations this gives a differential inequality
\begin{equation}   \Delta \vert \Riem \vert \geq - \vert \Riem\vert^{2}. \end{equation}
(We use the \lq\lq analysts convention'' for the sign of $\Delta$.)
Now apply this to our cone $Y$ and set $f = \vert \Riem\vert$.  Write $r$ for the radial function on the cone. Clearly $f$ is homogeneous of degree $-2$, so (in an obvious notation)
$$   \frac{\partial f}{\partial r} = -\frac{2}{r} f \ \ \ \ \frac{\partial^{2} f}{\partial r^{2}}= \frac{6}{r^{2}} f. $$
 Consider the restriction of $f$ to the cross-section $r=1$ and a point $p$ where $f$ attains its  maximal value $m$, say. If we write $\Delta_{\Sigma}$ for the Laplacian on the cross-section we have the usual formula
$$  \Delta= r^{-5} \frac{\partial}{\partial r} \left( r^{5} \frac{\partial}{\partial r}\right) + r^{-2} \Delta_{\Sigma}. $$
By the homogeneity of $f$ we have, at the point $p$,
$$ \Delta f =-4f+ \Delta_{\Sigma} f, $$
and $\Delta_{\Sigma} f\leq 0$ at $p$ by the maximum principle. So we deduce that $4m\leq m^{2}$ and if $m<4$ we must have $m=0$.

\subsection{The main argument}

The second foundation of our proof is the existence, due to Cheeger and Colding \cite{kn:CC}, of \lq\lq tangent cones at infinity'', under suitable hypotheses, as metric cones.  
It is convenient to first  state an alternative form of Theorem 2.
\begin{prop}
With the same hypotheses as in Theorem 2, 
for all $\epsilon>0$ there is an $A(\epsilon)$ such that if the ball
$B(x,r)\subset M$ does not carry homology and $E(x,r)\leq \epsilon A(\epsilon)^{-2}$ then for all $\rho\leq r$ we have $E(x,\rho)\leq \epsilon$.
\end{prop}
To see that this statement implies Theorem 2, let $B(y,r')\subset B(x,r)$ be a ball of the kind considered there. The $r/2$ ball $B(y,r/2)$ centred at $y$ lies in $B(x,r)$ and the corresponding normalised energy is at most $2^{2}$ times $E(x,r)$.  It suffices to take
$\delta(\epsilon)= \epsilon (2A(\epsilon))^{-2}$ and apply Proposition 6 to $B(y,r/2)$.

Now we start the proof of Proposition 6.  Notice that in proving the result we can obviously assume that $\epsilon$ is as small as we please and we will always suppose that $\epsilon$ is less than the constant $\theta$ of Proposition 5.  First fix $\epsilon, A$ and suppose that the statement is false for these parameters, so we have a point $x$ in a manifold $M$ and radii $\rho<r$ with
$E(x,\rho)\geq \epsilon, E(x,r)\leq \epsilon A^{-2}$. This implies that $\rho<A^{-1} r$. Choose the largest possible value of $\rho$, so $E(x, s)<\epsilon$ for $s\in (\rho,r]$ but $E(x,\rho)=\epsilon$. Still keeping $\epsilon$ fixed we suppose that we have  such violating examples for a sequence $A_{n}\rightarrow \infty$, and data $M_{n}, g_{n},  x_{n}, \rho_{n}, r_{n}$. Rescale the (Riemannian) metrics by the factor $\rho_{n}^{-1/2}$: without loss of generality there is a pointed Gromov-Hausdorff limit
$$     (M_{n}, \rho_{n}^{-1/2}g_{n}, x_{n}) \rightarrow (Z,O). $$

Here $Z$ may be singular. Let $R_{m}$ be a any sequence with $R_{m}\rightarrow \infty$. Then we can rescale the metric on $Z$ by factors $R_{m}^{-1}$ and  the general theory of \cite{kn:CC} (using the noncollapsing condition (4)) tells us that a subsequence converges to a \lq\lq tangent cone at infinity'' $Y$ which is a metric cone. The convergence is again in the sense of pointed Gromov-Hausdorff limits. 

Now go back to the smooth Riemannian manifolds $M_{n}$ and consider the rescalings $(R_{m} \rho_{n})^{-1/2} g_{n}$, indexed by $n,m$. If we choose a suitable function $n(m)$ which increases sufficiently rapidly then the corresponding sequence of Riemannian manifolds converges to the cone $Y$.  Call these based manifold $M_{m}, \tg_{m}, x_{m}$. We can also suppose $n(m)$ increases so rapidly that $ A_{n(m)}/R_{m}$ tends to infinity with $m$.
The choice of parameters means that
$ E(x_{m}, s, \tg_{m})< \epsilon <\theta $ if $s\leq A_{n}/R_{m}$ and $s>R_{m}^{-1}$. In particular this holds for $s=2$ and we are in  the position considered in Theorem 2. Thus we deduce that the tangent cone $Y$ is smooth (away from its vertex) and of the form $\bC^{3}/\Gamma$. It follows easily that $Z$ is smooth outside a compact set and  that the curvature of $Z$ satisfies a bound
\begin{equation} \vert \Riem\vert \leq \sigma r_{Z}^{-2}, \end{equation}
(outside a compact set) where $r_{Z}$ denotes the distance in $Z$ to the base point and the constant $\sigma$ can be made as small as we like by choosing $\epsilon$ small.

The further result that we need is that the curvature actually decays faster.
\begin{prop}
There is a $\gamma>3$ and $K>0$ such that outside a compact subset of $Z$.
$$ \vert \Riem \vert \leq K r_{Z}^{-\gamma}. $$
\end{prop}

In fact our proof will establish the result for any $\gamma< 4$. The  decay rate $r^{-4}$ arises as that of  the Green's function in real dimension $6$. We postpone the proof of Proposition 7, which is somewhat standard,  and move on to complete the proof of Theorem 3, assuming this.

Consider the general situation of a domain $U\subset M$ with smooth boundary
$\partial U$ and a K\"ahler metric on $M$. Suppose that 
\begin{itemize}
\item The cohomology class of $\omega$ in $H^{2}(\partial U)$ is zero.
\item  $U$ does not carry homology.
\item $H^{1}(\partial U)=0$.
\end{itemize}

 The first two conditions imply that the class of $\omega$ in $H^{2}(U)$ also vanishes, so we can write $\omega=da$ for a $1$-form $a$ on $U$. 
Now let $p$ be the invariant polynomial corresponding to the characteristic class $a_{3} c_{1}^{2}+ b_{3} c_{2}$, as discussed in the Introduction. Thus
$p(\Riem)$ is a closed $4$-form on $U$ and if the metric is K\"ahler-Einstein we have
$$  \int_{U} \vert \Riem\vert^{2}= \int_{U} p(\Riem)\wedge \omega. $$
 Then
\begin{equation}  \int_{U} p(\Riem) \wedge \omega = \int_{\partial U} p(\Riem) \wedge a. \end{equation}
Further, the third condition implies that if $\tilde{a}$ is {\it any} $1$-form on $\partial U$ with $d\tilde{a}=\omega\vert_{\partial U}$ then
\begin{equation}   \int_{\partial U}p(\Riem) \wedge a =  \int_{\partial U}p(\Riem) \wedge\tilde{a}.\end{equation}
Thus the integral of $p(\Riem)\wedge \omega$ over $U$ is determined by data on the boundary and if the metric is K\"ahler-Einstein this coincides with
the integral of $\vert \Riem\vert^{2}$ over $U$. In particular we get an inequality
\begin{equation} \int_{U} \vert \Riem \vert^{2} \leq c \max_{\partial U}\vert \Riem\vert^{2} \  \int_{\partial U}\vert a \vert, \end{equation}
for some fixed constant $c$.

To apply this, consider first the flat cone  $\bC^{3}/\Gamma$ and let   $\Sigma_{r}$ be the cross-section at radius $r$. Then the real cohomology of $\Sigma_{r}$ vanishes in dimensions $1$ and $2$ and we can obviously write
$$  \omega\vert_{\Sigma_{r}} = d a_{r} $$
with $\vert a_{r}\vert = O(r)$. In fact we can take $a_{r}$ to be one half the contraction of $\omega$ with the vector field $ r \frac{\partial }{\partial r}$, so $\vert a_{r}\vert = r$. Let $V_{\Gamma}\  r^{5}$ be the 5-volume of $\Sigma_{r}$.

Now choose $R_{0}$ so large  that $c K^{2} R_{0}^{6-2\gamma}V_{\Gamma} <\epsilon/2$.
By our discussion of the cone at infinity in $Z$ we can choose $R\geq R_{0}$ such that near the level set $r_{Z}=R$ in $Z$ there is a hypersurface
$\Sigma'$ with the property that the geometry of $\omega_{Z}$ restricted to $\Sigma'$ is very close to that of the cone metric restricted to $\Sigma_{R}$, in an obvious sense. In particular we can suppose that $\omega_{Z}\vert_{\Sigma'}= da'$ where
$$    \int_{\Sigma'} \vert a' \vert \leq \frac{3}{2} V_{\Gamma} R^{6}, $$
say. The bound on the curvature in Proposition 7  implies that
\begin{equation} c  \max_{\Sigma'} \vert \Riem \vert^{2}  \int_{\Sigma'} \vert a'  \vert \leq \frac{3\epsilon}{4}. \end{equation}

From now on $R$ is fixed. We go back to the manifolds $M_{n}$ with rescaled metrics $\rho_{n}^{-1/2} g_{n}$ converging to $Z$. We choose $n$ so large that there is hypersurface
$\Sigma''\subset M_{n}$ on which the geometry is close to that of $\Sigma'$. Then $\Sigma''$ is the boundary of an open set $U\subset M_{n}$ which does not carry homology by hypothesis. By taking $n$ large we can make the boundary term in (10) as close as we like to that estimated in (12),(13). So we can suppose that
$$  \int_{U} \vert \Riem\vert^{2} \leq \frac{5\epsilon}{6}, $$
say. But the  $U$ contains the unit ball centred at $x_{n}$ over which the integral of $\vert \Riem\vert^{2}$ is $\epsilon$ by construction. This gives the desired contradiction.

\subsection{Curvature decay}

Here we prove Proposition 7. The proof exploits the differential inequality $$\Delta \vert \Riem \vert \geq -\vert \Riem \vert^{2} \geq - \sigma r_{Z}^{-2} \vert \Riem \vert. $$
To explain the argument consider first a slightly different problem in which we work on a cone $Y$ (of real dimension 6) with radial function $r$ and we have a smooth positive function $f$ on the set $r\geq 1$ with $\Delta f \leq \sigma r^{-2} f$ and $f\leq \sigma r^{-2}$. We have
\begin{equation}   (\Delta +\frac{\sigma}{r^{2}}) r^{\lambda}= 
( \lambda(\lambda+4)+\theta) r^{\lambda-2}. \end{equation}
Let $\alpha=-2-\sqrt{4-\sigma}, \beta=-2+\sqrt{4-\sigma}$: the roots of the equation $\lambda(\lambda+4)+\theta=0$. We suppose $\sigma$ is small, so $\alpha$ is close to $-4$ and $\beta$ is close to $0$. Then
any linear combination
$$   g= A r^{\alpha} + B r^{\beta} $$
satisfies the equation $(\Delta+ \frac{\sigma}{r^{2}}) g =0$. Fix $R>1$ and let $m_{1}, m_{R}$ be the maximum values of $f$ on the cross-sections $r=1, r=R$ respectively. We choose constants $A,B$ so that $g(1)\geq m_{1}, g_{R}\geq m_{R}$.
If we solve for the case of equality we get
$$   A_{0}= \frac{ R^{\beta} m_{1}- m_{R}}{ R^{\beta}-R^{\alpha}}\ ,\ B_{0}= \frac{m_{R}- R^{\alpha} m_{1}}{R^{\beta} - R^{\alpha}}. $$
Since $R^{\beta}>R^{\alpha}$  it certainly suffices to take
$$  A= \frac{R^{\beta} m_{1}}{R^{\beta}-R^{\alpha}}\ ,\ B=\frac{m_{R}}{R^{\beta}-R^{\alpha}}. $$
Then set $u=f-g$ so that $u\leq 0$ when $r=1$ or $r=R$ and $(\Delta+\sigma r^{-2})u\geq 0$. We claim that $u\leq 0$ throughout the region $1\leq r\leq R$. For if we write $u= h r^{-2}$, calculation gives
\begin{equation}  (\Delta+\frac{\sigma}{r^{2}}) u = r^{-2} (\Delta h - 2 r^{-1} \frac{\partial h}{\partial r} - (4-\sigma) r^{-2} h) \end{equation} and (since $4-\sigma>0$) we see that there can be no interior positive maximum of $h$.

We see then that for $1< r<R$ we have
$$    f\leq A r^{\alpha}+ B r^{\beta} $$
where $A,B$ are given by the formulae above. Using the information $m_{1}\leq \theta, m_{R}\leq \sigma R^{-2}$ we get
$$  Ar^{\alpha} \leq c r^{\alpha} \ \ , B r^{\beta} \leq c R^{-\beta+2} r^{\beta},$$
for some fixed $c$. Taking $R$ very large compared to $r$ we get $f\leq 2c r^{\alpha}$, say.

We want to adapt this argument to the function $f=\vert \Riem\vert$ on $Z$.
Recall that $Z$ has base point $O$ and $r_{Z}$ is the distance to $O$. We write
$$  \cA(r_{1}, r_{2})= \{ z\in Z: r_{1}< r_{Z}(z)< r_{2} \}$$
with $r_{2}=\infty$ allowed. 

 There are several complications. One minor difficulty is that $f$ may not be smooth, but this handled by standard approximation arguments. The second is that the manifold is not exactly a cone, even at large distances, and the radius function $r_{Z}$ need not be smooth.
\begin{lem}
For any $\tau>0$ we can find an $R_{\tau}>0$ and a smooth function $\ur$ on the region
$\cA(R_{\tau},\infty)\subset Z$  such that 
\begin{itemize}
\item $\vert \frac{\ur}{r_{Z}}-1\vert \leq \tau, $
\item $ \vert \vert \nabla \ur \vert-1\vert \leq \tau$,
\item $ \vert \Delta \ur^{2} - 12 \vert \leq \tau$.
\end{itemize}
\end{lem}

We only sketch the proof. It is clear that we can choose $R_{\tau}$ so that find such a function on any annulus $\cA(2^{p} R_{\tau}, 2^{p+1} R_{\tau})$ for $p=1,2,\dots$. This just uses the convergence of the rescaled metric to the cone. If we extend these annuli slightly we get a sequence of overlapping annuli and a function defined on each. To construct $\ur$ we glue these together using cut-off functions. Notice that we only have to glue adjacent terms so that the gluing errors do not accumulate.

Now think of $\sigma$ and $\tau$ as fixed small numbers.( It will be clear from the discussion below that one could calculate appropriate values explicitly: for example $\sigma=\tau=1/100$ will do) We want to adapt the preceding argument to prove Proposition 7. There is no loss in supposing that in fact $R_{\tau}=1/2$ and that $\vert \Riem \vert \leq \sigma \vert \Riem \vert^{2}$ on $\ur\geq 1$. Let $\alpha, \beta$ be the roots as above and choose $\alpha',\beta'$ with $\alpha'$ slightly greater than $\alpha$ and $\beta'$slightly less than $\beta$. Then if $\tau$ is small we will have
\begin{equation}  (\Delta+\sigma \ur^{-2})\ur^{\alpha'}\ ,\  (\Delta+\sigma \ur^{-2})\ur^{\beta'}\leq 0. \end{equation}
We want to choose $A,B$ such that $g(\ur)=A\ur^{\alpha'}+B \ur^{\beta'}$ has
$\Delta+ \sigma \ur^{-2}\leq 0$ and $g(1)\geq m_{1}, g(R)\geq m_{R}$ where
$m_{1}, M_{R}$ are the maxima of $f$ on $\ur=1, \ur=R$ respectively. We take
$$   A= \frac{ m_{1}R^{\beta'}}{R^{\beta'}-R^{\alpha'}}\ ,\  B=\frac{m_{R}}{R^{\beta'}-R^{\alpha'}}.
$$
Then $A,B\geq 0$ so the differential inequality follows from (16). Now consider $u=f-g$ and argue as before to show that $u$ has no interior maximum. We have
$$  \Delta \ur^{-2}= \ur^{-4} ( 8 \vert \nabla \ur\vert^{2}-\Delta \ur^{2})\leq (-4+9\delta )\ur^{-4}, $$
and the same argument goes through.

\section{Proof of Theorem 1}
We begin with a standard result.
\begin{prop}
Given $\kappa>0$ there is a $\chi>0$ such that if $B$ is any unit ball in a K\"ahler-Einstein manifold (of real dimension $6$) and
\begin{enumerate}\item  $\vert \Riem \vert =1$ at the center of $B$;
 \item $\vert \Riem \vert \leq 4$ throughout $B$;
 \item $\Vol(B)\geq \kappa$
 \end{enumerate}
 then $\int_{B} \vert \Riem \vert^{2} \geq \chi$.
 \end{prop}
  One way to prove this is to apply the  Moser iteration technique to $\vert\Riem \vert$, using the differential inequality (8) and the fact that in this situation the Sobolev constant is bounded. Another method is to use elliptic estimates in harmonic coordinates.

 \begin{prop}
 Suppose $M$ is a K\"ahler-Einstein manifold as considered in Theorem 1. There is an $\epsilon>0$ such that$U\subset M$ is any domain such that  the normalised energy of any ball contained in $U$ is less than $\epsilon$ then $\vert \Riem \vert \leq 4 d^{-2}$, where $d$ denotes the distance to the boundary of $U$.
\end{prop}

To see this we let $S$ be the maximum value of $d^{-2} \vert Riem \vert$ over $U$ and suppose that this is attained at $p$. If $S>4$ then $d(q)\geq d(p)/2$ for any point $q$ in the ball of radius $d(p)S^{-1/2}$ centred at $p$. Rescale this ball to unit size and we are in the situation considered in the preceding proposition. If we take $\epsilon$ to be the constant $\chi$ appearing there we get a contradiction, so in fact $S\leq 4$. 

We can now prove Proposition 1. Take $\epsilon$ as above and let $\delta(\epsilon)$ be the value given by Theorem 3. Suppose that $B(x,r)$ is a ball of normalised energy less than $\delta$ and let $U$ be the half-sized ball. Then Theorem 3 tells us that the normalised energy of any ball in $U$ is at most $\epsilon$ and we can apply the result above to see that $\vert \Riem \vert \leq 4.6^{2} r^{-2}$ in $B(x,r/3)$.

\

Of course some of the constants appearing in the statements above (such as the value $4$ in Proposition 8) are rather arbitrary. The essential point is that there is a definite threshold value, so that if the normalised energy is below this threshold, on a ball which does not carry homology, then we get complete control of the metric on interior regions.

\

Now we prove  Theorem 1. Given $r$ and a K\"ahler-Einstein metric $(M,g)$ as considered there we pick a maximal collection of point $x_{a}, a\in I $ in $M$ such that the distance between any pair is at least $r$. Then the $r$-balls
$B(x_{a}, r)$ cover $M$. Consider a ball $B(x_{a}, 12 r)$. If this ball does not carry homology  and its normalised energy is less than $\delta$ then by the result above we have
$\vert \Riem \vert \leq 1$ on $B(x_{a}, 2r)$. Thus no point in $Z_{r}$ can lie in   $B(x_{a}, r)$. Thus $Z_{r}$ is covered by balls $B(x_{a}, r)$ where
either $E(x_{a}, 12 r) >\delta$ or $B(x_{a}, 12 r)$ carries homology. Let
$I'\subset I$ denote the indices of the first kind and $I''\subset I$ those of the second.

Suppose that $N$ balls of the form $B(x_{a}, 12 r)$ have a non-empty common intersection. Let $q$ be a point in the intersection, so the $N$ centres $x_{a}$ all lie in the $12r$ ball centred at $q$. By construction the balls $B(x_{a}, r/2)$ are disjoint and have volume at least $(\kappa/2^{6}) r^{6}$.
Since the volume of $B(q,12.5 r)$ is bounded above by a fixed multiple of $r^{6}$ this gives a fixed bound on $N$, independent of $r$. 
Thus
$$  \sum_{a\in I'}\int_{B(x_{a}, 12 r)} \vert \Riem\vert^{2} \leq N \int_{M} \vert \Riem \vert^{2}. $$
On the other hand, by definition,  $E(x_{a}, 12 r) \geq \delta$ for $a\in I'$ so we see that the number of elements of $I'$ is at most
$$  \vert I'\vert \leq \frac{N}{\delta (12 r)^{2}} \int_{M} \vert \Riem\vert^{2}= C_{1} r^{-2},$$
say.

By a similar argument there is a fixed upper bound $N'$ on the number of balls $B(x_{a}, 12 r)$ which can meet any given one. It follows that the number of these balls which carry homology is bounded by $N'$ times the second Betti number of $M$. So $\vert I''\vert$ is bounded by a fixed number. Then $$ \Vol(Z_{r}) \leq {\rm constant}\  r^{6} (\vert I' \vert + \vert I''\vert)\leq C_{1} r^{4}+ C_{2} r^{6}. $$

\

Now we turn to Theorem 3. Recall that this states that we can find a {\it connected } open subset $\Omega'\subset M\setminus Z_{r}$ so that the volume of the complement of $\Omega'$ is bounded by a multiple of $r^{18/5}$. 

To prove this we recall that in our situation there is a bound on the isomperimetric constant, due to Croke \cite{kn:Cr}. If $H\subset M$ is a rectifiable hypersurface dividing $M$ into two components $M_{1}, M_{2}$ with $\Vol(M_{1})\leq \Vol(M_{2})$ then
\begin{equation}    \Vol(M_{1})^{1/6} \leq k \Vol(H)^{1/5}, \end{equation}
for a fixed constant $k$. By the construction in the proof of Theorem 1 above the set $Z_{r}$ is contained in $W$ which is a union of $P$ balls of radius $r$ with $P\leq C (r^{-2}+1)$. By the Bishop comparison theorem the $5$-volume of the boundary of one of these balls is bounded by $C r^{5}$. The boundary $\partial W$ of $W$ is a rectifiable set of $5$-volume at most the sum of the boundaries of the balls, thus $\Vol(\partial W) \leq C( r^{3}+ r^{5})$. We can normalise so that the  the volume of $M$ is $1$ and then, without loss of generality suppose that $r$ is so small that $\Vol(W)\leq 1/10$ and $k^{5}\Vol(\partial W)\leq 9/10$.
Let $\Omega_{i}$ be the connected components of $M\setminus W$. If a component
$\Omega_{1}$ has volume greater than $1/2$ (i.e. one half the volume of $M$) then by (14) its complement has volume less than $k^{6} \Vol(\partial W)^{6/5}= = k^{6} C^{6/5} r^{18/5}$ and we can take $\Omega'=\Omega_{1}$. So suppose that all components $\Omega_{i}$ have volume less than $1/2$. Then it is clear that
   $$\sum_{i} \Vol(\partial \Omega_{i}) = \Vol(\partial W)$$
   while $$\sum_{i} \Vol(\Omega_{i})\geq 9/10.$$
   The second equation implies that $\sum_{i} \Vol(\Omega_{i})^{5/6}\geq 9/10$ and then (17) gives
$$  \sum_{i} \Vol(\partial \Omega_{i})\geq 9 k^{-5}/10, $$
so
$$ \Vol(\partial W) \geq 9k^{-5}/10, $$
contrary to our assumptions.

\section{Discussion}
\begin{enumerate}
\item Proposition 5 and the ensuing arguments in 2.2 above are closely related to a result of Cheeger, Colding and Tian (\cite{kn:CCT}, Theorem 9.26). Let $X$ be a complete, noncompact Ricci-flat K\"ahler  manifold of complex dimension $n$ with base point $p$. Suppose that ${\rm Vol}(B(p,R)\geq \kappa R^{2n}$ and 
$$ R^{2-n} \int_{B(p,R)} \vert \Riem\vert^{2} \rightarrow 0 $$ as $R\rightarrow \infty$.  Tian conjectured in \cite{kn:T} that  in this situation $X$ is an ALE manifold, with tangent cone at infinity of the form  $\bC^{n}/\Gamma$. When $n=3$ this conjecture was established in the result quoted above. Our statements are a little different since we establish a definite \lq\lq small asymptotic energy threshold'' which implies that the manifold is  ALE. 

\item In our situation the tangent cone $Y$ of $Z$ at infinity is unique.
In general, positive Ricci curvature does not imply uniqueness of tangent cones, even when $\vert \Riem \vert \leq c r^{-2}$.  See the discussion in \cite{kn:CC1} of examples, including an unpublished example due to Perelman. 

\item Theorem 3 becomes  false if we omit the condition that the ball does not carry homology. To see this one can consider for example the quotient
$\underline{M}= T^{6}/\Gamma$ of a complex torus by a group $\Gamma$ of order $3$, acting with isolated fixed points. Then $\underline{M}$ has a  resolution
$M$ with $c_{1}(M)=0$.  According to Joyce \cite{kn:Joyce} there is a family of ALE metrics on the resolution of $\bC^{3}/\Gamma$, parametrised by the K\"ahler class. For suitable K\"ahler classes on $M$, the Calabi-Yau metric is approximately given by gluing rescaled versions of  these ALE metrics to the flat metric on
$\underline{N}$, just as in the familiar picture of the Kummer construction for K3 surfaces.
For any $\delta>0$ we can use this  scaling to find a Ricci-flat metric on $M$, a unit ball $B(x,1)\subset M$ on which the normalised energy is less than $\delta$ and an interior ball $B(x,\rho)$ on which the normalised energy exceeds $c$, for some fixed $c$. The arguments in our forthcoming paper will  will show that this is essentially the only way in which approximate monotonicity can fail; see also \cite{kn:T1}. 
\item It is interesting to ask if similar results to those proved above can be established for  constant scalar curvature and extremal K\"ahler metrics. For this it might be sensible to assume a bound on the Sobolev constant. Perhaps some of the  techniques  used in \cite{kn:TV}, \cite{kn:CW} can be applied to this problem. 
\item Another question is whether a result like Theorem 3 holds for general 6-dimensional Einstein metrics (with nonnegative Ricci curvature and a non-collapsing condition). It might be that a different topological side condition is appropriate.

\item The function $\delta(\epsilon)$ in Theorem 3 depends only on the collapsing constant $\kappa$. It would be interesting to determine the function effectively, but our method cannot do this. One suspects that, in reality, it may be possible to take $\delta(\epsilon)$ not much smaller than $\epsilon$, and that the constant $C$ in Theorem 1 can (in reality) be taken not too large. Some results on the numerical analysis of K\"ahler-Einstein, and more generally extremal, metrics seem to give evidence for this suspicion \cite{kn:BD}, \cite{kn:DKLR} but a theoretical derivation of realistic estimates seems a long way off.  
\end{enumerate}


\end{document}